\documentclass[11pt]{article}
\usepackage{amsmath}
\usepackage{amssymb}
\usepackage{amsopn}
\usepackage{amscd}
\usepackage{epsfig}
\usepackage{amsfonts}
\usepackage{latexsym}
\usepackage{graphicx}
\usepackage{color}

\newtheorem{theorem}{Theorem}[section]

\newtheorem{lemma}{Lemma}[section]
\newtheorem{proposition}{Proposition}[section]

\newcommand{\ravphi}{{\varphi_{\textsc{\tiny R}}}}
\newcommand{\ate}{{\eta_{\textsc{\tiny R}}}}
\newcommand{\f}{w_{\rm F}}
\newcommand{\xf}{x_{\rm F}}
\newcommand{\Lg}{{\cal L}}

\definecolor{trp}{rgb}{1,1,1}

\definecolor{red}{rgb}{1,0,.2}

\definecolor{blue}{rgb}{0,0,1}

\definecolor{rgrey}{rgb}{.8,0.4,.4}  
\definecolor{grey}{rgb}{.13,.13,.13}  

\definecolor{green}{rgb}{0.0,0.4,0.2}

\begin{document}

\begin{center}

\medskip

{\huge{\bf On the conjugacy class of the Fibonacci \\[.3cm]dynamical system}}

\bigskip

\bigskip

{\large{Michel Dekking (Delft University of Technology)\\ and\\ Mike Keane (Delft University of Technology and University of Leiden)}}

\medskip

{\large{{Version: August 16, 2016}}}  

\bigskip

\bigskip

\end{center}



\section{Introduction}\label{sec:intro}

We study the Fibonacci substitution $\varphi$ given by
$$\varphi:\quad 0\rightarrow\,01,\;1\rightarrow 0.$$
The infinite Fibonacci word $\f$ is the unique one-sided sequence (to the right) which is a fixed point of $\varphi$:
$$\f=0100101001\dots.$$
We also consider one of the two two-sided fixed points $\xf$ of $\varphi^2$:
$$\xf=\dots01001001\!\cdot\!0100101001\dots.$$
The dynamical system generated by taking the orbit closure of $\xf$ under the shift map $\sigma$ is denoted by  $(X_\varphi,\sigma)$.

The question we will be concerned with is: what are the substitutions $\eta$ which generate a symbolical dynamical system topologically isomorphic to the Fibonacci dynamical system? Here topologically isomorphic means that there exists a homeomorphism $\psi: X_\varphi\rightarrow X_\eta$,  such that $\psi\sigma=\sigma\psi$, where we denote the shift on $X_\eta$ also by $\sigma$. In this case $(X_\eta, \sigma)$ is said to be conjugate to
$(X_\varphi,\sigma)$.

This question has been completely answered for the case of constant length substitutions in the paper \cite{CDK}.  It is remarkable that there are only finitely many injective primitive substitutions of length $L$ which generate a system conjugate to a given substitution of length $L$. Here a substitution $\alpha$ is called  \emph{injective} if $\alpha(a)\ne \alpha(b)$ for all letters $a$ and $b$ from the alphabet with $a\ne b$. When we extend to the class of all substitutions, replacing $L$ by the Perron-Frobenius eigenvalue of the incidence matrix of the substitution, then the conjugacy class can be infinite in general. See \cite{Dekking-TCS} for the case of the Thue-Morse substitution. In the present paper we will prove that there are infinitely many injective primitive substitutions with  Perron-Frobenius eigenvalue $\Phi=(1+\sqrt{5})/2$ which generate a system conjugate to the Fibonacci system---see Theorem~\ref{th:inf}.

In the non-constant length case some new phenomena appear. If one has an injective substitution $\alpha$ of constant length $L$, then all its powers $\alpha^n$ will also be injective. This is no longer true in the general case. For example, consider the injective substitution $\zeta$ on the alphabet $\{1,2,3,4,5\}$ given by
 $$\zeta: \qquad 1\rightarrow 12,\;
 2\rightarrow 3,\;
 3\rightarrow 45,\;
 4\rightarrow 1,\;
 5\rightarrow 23.$$
 An application of Theorem~\ref{th:Nblock} followed by a partition reshaping (see Section~\ref{sec:reshaping}) shows that the system $(X_\zeta,\sigma)$ is conjugate to the Fibonacci system.
  However, the square of $\zeta$ is given by
 $$\zeta^2: \qquad 1\rightarrow 123,\;
 2\rightarrow 45,\;
 3\rightarrow 123,\;
 4\rightarrow 12,\;
 5\rightarrow 345, $$
 which is \emph{not} injective. To deal with this undesirable phenomenon we introduce the following notion. A substitution $\alpha$ is called a \emph{full rank} substitution if its incidence matrix has full rank (non-zero determinant). This is a strengthening of injectivity, because obviously a substitution which is not injective can not have full rank. Moreover, if the substitution $\alpha$ has full rank, then all its powers $\alpha^n$ will also have full rank, and thus will be injective.

 Another phenomenon, which does not exist in the constant length case, is that non-primitive substitutions $\zeta$ may generate uniquely defined minimal systems
 conjugate to a given system. For example, consider the injective substitution $\zeta$ on the alphabet $\{1,2,3,4\}$ given by
 $$\zeta:\qquad 1\rightarrow 12,\quad
 2\rightarrow 31,\quad
 3\rightarrow 4,\quad
 4\rightarrow 3. $$
With the partition reshaping technique from Section~\ref{sec:reshaping} one can show that the system $(X_\zeta,\sigma)$ is conjugate to the Fibonacci system (ignoring the system on two points generated by $\zeta$). In the remainder of this paper we concentrate on primitive substitutions.

The structure of the paper is as follows. In Section~\ref{sec:Nblock} we show that all systems in the conjugacy class of the Fibonacci substitution can be obtained by letter-to-letter projections of the systems generated by so-called $N$-block substitutions. In Section~\ref{sec:C3} we give a very general characterization of symbolical dynamical systems in the Fibonacci conjugacy class, in the spirit of a similar result on the Toeplitz dynamical system in \cite{CKL08}. In Section~\ref{sec:reshaping} we introduce a tool which admits to turn non-injective substitutions into injective substitutions. This is used in Section~\ref{sec:C1} to show that the Fibonacci class has infinitely many primitive injective substitutions as members. In Section~\ref{sec:two} we quickly analyse the case of a 2-symbol alphabet. Sections \ref{sec:equi} and \ref{sec:mat} give properties of maximal equicontinuous factors and incidence matrices, which are used to analyse the 3-symbol case in  Section \ref{sec:C2}. In the final Section \ref{sec:L2L} we show that the system obtained by doubling the 0's in the infinite Fibonacci word is conjugate to the Fibonacci dynamical system, but can not be generated by a substitution.

\section{$N$-block systems and $N$-block substitutions}\label{sec:Nblock}

For any $N$ the $N$-block substitution $\hat{\theta}_N$ of a  substitution $\theta$ is defined on an alphabet of $p_\theta(N)$ symbols, where $p_\theta(\cdot)$ is the complexity function of the language $\Lg_\theta$ of $\theta$ (cf.\ \cite[p.~95]{Queff}). What is \emph{not} in \cite{Queff}, is that this $N$-block substitution generates the $N$-block presentation of the system $(X_\theta,\sigma)$.
 
 We denote the letters of the alphabet of the $N$-block presentation by $[a_1a_2\dots a_N]$, where $a_1a_2\dots a_N$ is an element from $\Lg_\theta^N$, the set of words of length $N$ in the language of $\theta$. The $N$-block presentation $(X^{[N]}_\theta,\sigma)$ emerges by applying an sliding block code $\Psi$ to the sequences of $X_\theta$, so $\Psi$ is the map \\[-.3cm]
 $$\Psi(a_1a_2\dots a_N)=[a_1a_2\dots a_N].$$
 We denote by $\psi$ the induced map from $X_\theta$ to $X^{[N]}_\theta$:
 $$\psi(x)=\dots\Psi(x_{-N},\dots,x_{-1})\Psi(x_{-N+1},\dots,x_{0})\dots.$$
 It is easy to see that $\psi$ is a conjugacy, where the inverse is $\pi_0$ induced by the 1-block map (also denoted $\pi_0$) given by  $\pi_0([a_1a_2\dots a_N])=a_1$.

 The $N$-block substitution  $\hat{\theta}_N$ is defined by requiring that for each word $a_1a_2\dots a_N$ the length of $\hat{\theta}_N([a_1a_2\dots a_N])$  is equal to the length  $L_1$ of $\theta(a_1)$, and the letters of $\hat{\theta}_N([a_1a_2\dots a_N])$ are the $\Psi$-codings of the first $L_1$ consecutive $N$-blocks in $\theta(a_1a_2\dots a_N)$.
 

 \begin{theorem}\label{th:Nblock} Let $\hat{\theta}_N$  be the $N$-block substitution of a primitive substitution $\theta$. Let $(X^{[N]}_\theta,\sigma)$ be the $N$-block presentation of the system $(X_\theta,\sigma)$. Then the set $X^{[N]}_\theta$ equals $X_{\hat{\theta}_N}$.
\end{theorem}

%
%

{\em Proof:} Let $x$ be a fixed point of $\theta$, and let $y=\psi(x)$, where $\psi$ is the $N$-block conjugacy, with inverse $\pi_0$.  The key equation is $\pi_0\,\hat{\theta}_N=\theta\,\pi_0$.
 This implies\\[-.3cm]
 $$\pi_0\,\hat{\theta}_N(y)=\theta\,\pi_0(y)=\theta\,\pi_0(\psi(x))=\theta(x)=x.$$

Applying $\psi$ on both sides gives $\hat{\theta}_N(y)=\psi(x)=y$, i.e., $y$ is a fixed point of $\hat{\theta}_N$. But then $X^{[N]}_\theta=X_{\hat{\theta}_N}$, by minimality of $X^{[N]}_\theta$. $\Box$

\medskip


It is well known (see, e.g., \cite[p.~105]{Queff}) that  $p_\varphi(N)=N+1$, so for the Fibonacci substitution $\varphi$ the $N$-block substitution $\hat{\varphi}_N$ is a substitution on an alphabet of $N+1$ symbols.

We describe how one obtains $\hat{\varphi}_2$.
We have $\Lg_\varphi^2=\{00, 01,  10\}$. Since 00 and 01 start with 0, and 10 with 1, we obtain
$$\hat{\varphi}_2:\quad [00]\mapsto [01][10],\;[01]\mapsto [01][10],\; [10]\mapsto [00],$$
reading off the consecutive 2-blocks from $\varphi(00)=0101,\, \varphi(01)=010$ and $\varphi(10)=001$.
It is  useful to recode the alphabet $\{[00],[01],[10]\}$ to the standard alphabet $\{1,2,3\}$.
We do this in the order  in which they appear for the first time in the infinite Fibonacci word $\f$--- we call this the \emph{canonical coding}, and will use the same principle for all $N$. For $N=2$ this gives
  $[01]\rightarrow 1,\; [10]\rightarrow 2,\; [00]\rightarrow 3$. Still using the notation $\hat{\varphi}_2$ for the 
 substitution on this new alphabet, we obtain
$$\hat{\varphi}_2(1)=12 \quad \hat{\varphi}_2(2)=3, \quad \hat{\varphi}_2(3)=12.$$
In this way the substitution is in standard form (cf.~\cite{CDK} and \cite{Dekking-2016}).

 \section{The Fibonacci conjugacy class}\label{sec:C3}


%
%

%
%
%
%

Let $F_n$ for $n=1,2,\dots$ be the Fibonacci numbers $$F_1=1,\, F_2=1,\, F_3=2,\, F_4=3,\, F_5=5, \dots.$$

 \begin{theorem} Let $(Y,\sigma)$ be any subshift. Then $(Y,\sigma)$ is topologically conjugate to the Fibonacci system $(X_\varphi,\sigma)$ if and only if there exist $n\ge 3 $ and two  words $B_0$ and $B_1$ of length $F_n$ and $F_{n-1}$,  such that any  $y$ from $Y$ is a concatenation of $B_0$ and $B_1$, and moreover, if\, $\cdots B_{x_{-1}} B_{x_0} B_{x_1}\cdots B_{x_k}\cdots$ is such a concatenation, then $x=(x_k)$ is a sequence from the Fibonacci system.
\end{theorem}

\noindent \emph{Proof:}  First let us suppose that $(Y,\sigma)$ is topologically isomorphic to the Fibonacci system. By the Curtis-Hedlund-Lyndon theorem, there exists an integer $N$ such that $Y$ is obtained by a letter-to-letter projection $\pi$ from the $N$-block presentation $(X^{[N]}_\varphi, \sigma)$ of the Fibonacci system.
Now if $B_0$ and $B_1$ are two decomposition blocks of sequences from
$X^{[N]}_\varphi$ of length $F_n$ and $F_{n-1}$, then $\pi(B_0)$ and $\pi(B_1)$ are decomposition blocks of sequences from $Y$ with lengths $F_n$ and $F_{n-1}$, again satisfying the concatenation property.
 So it suffices to prove the result for $X^{[N]}_\varphi$. Note that we may suppose that the integers $N$ pass through an infinite subsequence; we will use $N=F_n$,where $n=3,4,\dots$. Useful to us are the   \emph{singular words} $w_n$ introduced in \cite{WenWen}. The  $w_n$ are  the unique words of length $F_{n+1}$ having a different Parikh vector from all the other words of length $F_{n+1}$ from the language of $\varphi$. Here $w_1=1, w_2=00$, $w_3=101$, and for $n\ge4$
 $$w_n=w_{n-2}w_{n-3}w_{n-2}.$$
  The set of return words of $w_n$ has only two elements which are $u_n=w_nw_{n+1}$ and $v_n=w_nw_{n-1}$ (see page 108 in \cite{HuangWen}).
 The lengths of these words are $|u_n|=F_{n+3}$ and  $|v_n|=F_{n+2}$. Let $w_n^-$ be $w_n$ with the last letter deleted.
 Define for $n\ge5$
 $$B_0=\Psi(u_{n-3}w_{n-3}^-), \quad B_1=\Psi(v_{n-3}w_{n-3}^-),$$
 where $\Psi$ is the $N$-block code from $\Lg_\varphi^N$ to $\Lg_{\varphi^{[N]}}$, with $N=F_{n-2}$.
  Then these blocks have the right lengths, and by   Theorem 2.11 in \cite{HuangWen}, the two return words partition the infinite Fibonacci word $\f$ according to the infinite Fibonacci word---except for a prefix $r_{n,0}$:
 $$\f=r_{n,0}u_nv_nu_nu_nv_nu_n\dots.$$
 By minimality this property carries over to all two-sided sequences in the Fibonacci dynamical system.

 \smallskip

For the converse, let $Y$ be a Fibonacci concatenation system as above. Let $C_0=\varphi^{n-2}(0)$ and $C_1=\varphi^{n-2}(1)$. We define a map $g$ from  $(Y,\sigma)$ to  a subshift of $\{0,1\}^{\mathbb{Z}}$ by
$$g:\quad \cdots B_{x_{-1}} B_{x_0} B_{x_1}\cdots B_{x_k}\cdots\; \mapsto \; \cdots C_{x_{-1}} C_{x_0} C_{x_1}\cdots C_{x_k}\cdots,$$
respecting the position of the $0^{\rm th}$ coordinate. Since $|C_0|=|B_0|$ and $|C_1|=|B_1|$, $g$ commutes with the shift. Also, $g$ is obviously continuous.
Moreover, since for any sequence $x$ in the Fibonacci system $\varphi^{n-2}(x)$ is again a sequence in the Fibonacci system, $g(Y)\subseteq X_\varphi$.
So, by minimality, $(X_\varphi,\sigma)$  is a factor of  $(Y,\sigma)$. Since $g$ is invertible, with continuous inverse, $(Y,\sigma)$ is in the conjugacy class of the Fibonacci system. \hfill $\Box$

\medskip

\noindent {\bf Example}\; The case $(F_n,F_{n-1})=(13,8)$. Then $n=7$, so we have to consider the singular word $w_4=00100$ of length 5.

\noindent The set of $5$-blocks is
$\{01001,\,10010,\,00101,\,01010,\,10100,\,00100\}.$\\
These will be coded by the canonical coding $\Psi$ to the standard alphabet $\{1,2,3,4,5,6\}$.
Note that $\Psi(w_4)=6$.
Further, $w_3=101$ and $w_5=10100101$.  So $u_n=0010010100101$ and $v_n=00100101$. Applying $\Psi$ gives the two decomposition
blocks  $B_0 = 6123451234512$ and $B_1 = 61234512$.


\section{Reshaping substitutions}\label{sec:reshaping}

 We call a language preserving transformation of a substitution a reshaping. An example is the prefix-suffix change used in \cite{Dekking-TCS}.
Here we consider a variation which we call a \emph{partition reshaping}. 

We give an example of this technique.
Take the $N$-block representation of the Fibonacci system for  $N=4$. All five 4-blocks occur consecutively at the beginning of the  Fibonacci word $\f$ as
$\{0100,\,1001,\,0010,\, 0101,\, 1010\}.$
 The canonical coding  to $\{1,2,3,4,5\}$ gives the 4-block substitution $\hat{\varphi}_4$:
$$\hat{\varphi}_4:\qquad 1\rightarrow 12,\;
 2\rightarrow 3,\;
 3\rightarrow 45, \;
 4\rightarrow 12,\;
 5\rightarrow 3.$$

\noindent Its square  is equal to
$$\hat{\varphi}_4^2: \qquad 1\rightarrow 123,\;
 2\rightarrow 45,\;
 3\rightarrow 123,\;
 4\rightarrow 123,\;
 5\rightarrow 45. $$
 Since the two blocks $B_0=123$ and $B_1=45$  have no letters in common
 this permits to do a partition reshaping. Symbolically this can be represented by
 
 
\begin{table}[h!]
  \centering
  \caption{\small Partition reshaping.}
  \label{tab:table1}
  \begin{tabular}{ccccccccc}\\[.005cm]
  1            & \;       & 2            & 3           &     & \qquad\qquad     4            &  \;    & 5            &  \\
  $\downarrow$ & \;       & $\downarrow$ & $\downarrow$&     & \qquad\qquad     $\downarrow$ &  \;    & $\downarrow$ &  \\
   1           &  2       & 3            & 4           &  5  & \qquad\qquad     1            &  2     &  3           &  \\
   1           & \,\; 2 $\|$& 3            &\,\; 4 $\|$    &  5  & \qquad\qquad    \,\; 1  $\|$    &  2     & \,\; 3  $\|$   &
  \end{tabular}
\end{table}

Here the third line gives the images $\hat{\varphi}_4(B_0)=\hat{\varphi}_4(123)=12345$ and $\hat{\varphi}_4(B_1)=\hat{\varphi}_4(45)=123$; the fourth line gives a \emph{another} partition of these two words in three, respectively two subwords from which the new substitution $\eta$ can be read of:
$$\eta: \qquad 1\rightarrow 12,\;
 2\rightarrow 34,\;
 3\rightarrow 5 ,\;
 4\rightarrow 1,\;
 5\rightarrow 23. $$
What we gain is that the partition reshaped substitution $\eta$ generates the same language as $\hat{\varphi}_4$, but that $\eta$ is injective---it is even of full rank.

\section{The Fibonacci class has infinite cardinality}\label{sec:C1}

\smallskip

   \begin{theorem}\label{th:inf} There are infinitely many primitive injective substitutions with Perron-Frobenius eigenvalue the golden mean that generate dynamical systems topologically isomorphic to the Fibonacci system.
   \end{theorem}

 \smallskip

We will explicitly construct infinitely many primitive injective substitutions whose systems are topologically conjugate to the Fibonacci system. The topological conjugacy will follow from the fact that the systems are $N$-block codings of the Fibonacci system, where $N$ will run through the numbers $F_n-1$.
As an introduction we look at $n=5$, i.e., we consider the blocks of length $N=F_5-1=4$.
With the canonical coding of the $N$-blocks we obtain the 4-block substitution $\hat{\varphi}_4$---see Section~\ref{sec:reshaping}:
$$\hat{\varphi}_4:\qquad 1\rightarrow 12,\,
 2\rightarrow 3,\,
 3\rightarrow 45, \,
 4\rightarrow 12,\,\
 5\rightarrow 3.$$

\noindent An \emph{interval} $I$  starting with $a\in A$ is  a word of length $L$ of the form $$I=a,a+1,...,a+L-1.$$

\noindent Note that $\hat{\varphi}_4(123)=12345$, and   $\hat{\varphi}_4(45)=123$, and these four words are intervals.

This is a property that holds in general. First we need the fact that the first $F_n$ words of length $F_n-1$ in the
fixed point of $\varphi$ are all different. This result is given by Theorem 2.8 in \cite{Chuan-Ho}. We code these $N+1$ words by the canonical coding
to the letters $1,2,\dots,F_n$. We then have
\begin{equation}\label{eq:Fib}\hat{\varphi}_N(12...F_{n-1})=12\dots F_{n}, \qquad \hat{\varphi}_N(F_{n-1}\!+1,\dots F_n)=12\dots F_{n-1}.\end{equation}
This can be seen by noting that $\pi_0 \hat{\varphi}_N^n=\varphi^n \pi_0,$ for all $n$, and that the fixed point of $\varphi$ starts with $\varphi^{n-2}(0)\varphi^{n-3}(0)$.

\medskip

We continue  for  $n \ge 5$ with the construction of a substitution $\eta=\eta_n$ which is a partition reshaping of $\hat{\varphi}_N$.
The $F_n$ letters in the alphabet $A^{[N]}$ are divided in three species, S, M and L (for Small, Medium and Large).
$${\rm S}:={1,...,F_{n-3}}, \quad {\rm M}:={F_{n-3}+1,...,F_{n-1}},\quad  {\rm L}:={F_{n-1}\!+1,...,F_n}.$$
Note that ${\rm Card}\, {\rm M}=F_{n-1}-F_{n-3}=F_{n-2}=F_{n}-F_{n-1}={\rm Card}\, {\rm L}.$

An important role is played by $a_{ \rm M}:=F_{n-3}+1$, the smallest letter in M, and $a_{ \rm L}:=F_{n-1}+1$, the smallest letter in L.

For the letters in M (except for $a_{ \rm M})$ the rules are very simple:
$$\eta(a)= a+F_{n-2}$$ (i.e., a single letter obtained by addition of the two integers).
The first letter in M has the rule  $$\eta(a_{ \rm M})=\eta(F_{n-3}\!+1)= F_{n-1}, F_{n-1}\!+1= F_{n-1},a_{ \rm L} .$$
The images of the letters in L are intervals of length 1 or 2, obtained from a partition of the word $12\dots F_{n-1}$.
Their lengths are coming from $\varphi^{n-4}(0)$, rotated once (put the 0 in front at the back). This word is denoted $\rho(\varphi^{n-4 }(0))$.
The choice of this word is somewhat arbitrary, other choices would work. The properties of $v:=\rho(\varphi^{n-4}(0))$ which matter to us are

(V1) $\ell:=|v|=F_{n-2}$.

(V2) $v_1=1$, $v_\ell=0$.

(V3) $v$ does not contain any 11.

\noindent Now the images of the letters in L are determined by $v$ according to the following rule:  $|\eta(a_{ \rm L}+k-1)|=2-v_k$,   for all $k=1,\dots,F_{n-2}$.  Note that this implies in particular that for all $n\ge 5$ one has by property (V2)
$$\eta(a_{\rm L})=\eta(F_{n-1}\!+1)= 1, \;\qquad  \eta(F_n)= F_{n-1}-1,F_{n-1}.$$

The images of the letters in S are then obtained by choosing the lengths of the $\eta(a)$ in such a way that the largest common refinement of the induced partitions of the images of S and L is the singleton partition.

\medskip

\noindent{\bf Example}  The case $n=7$, so $ F_n=13$, $ F_{n-1}=8$, and $ F_{n-2}=5$.

\noindent Then ${\rm S}=\{1,2,3\},\, {\rm M}=\{4,5,6,7,8\},\, {\rm L}=\{9,10,11,12,13\}.$

\noindent Rules for M: \quad $4\rightarrow 89,\;5\rightarrow 10, \;6\rightarrow11, \;7\rightarrow 12, \;8\rightarrow13.$
Now $$\varphi^3(0)=01001\; \Rightarrow\; v= 10010\; \Rightarrow\; {\rm the\, partition\, is}\; 1|23|45|6|78.$$ This partition gives the following rules for L: $$9\rightarrow1,\; 10\rightarrow23,\; 11\rightarrow45,\; 12\rightarrow6,\; 13\rightarrow 78.$$
The induced partition for the images of the letters in S is $|12|34|567|8$, yielding rules
$$1\rightarrow12,\; 2\rightarrow34,\; 3\rightarrow567.$$

\noindent In summary we obtain the substitution $\eta=\eta_7$ given by : \vspace*{-0.1cm}
 \begin{align*}
{\rm S}:
\begin{cases}
1& \rightarrow 1,2\\
2& \rightarrow 3,4\\
3& \rightarrow 5,6,7
\end{cases}
\qquad {\rm M}:
\begin{cases}
4& \rightarrow 8,9\\
5& \rightarrow 10\\
6& \rightarrow 11\\
7& \rightarrow 12\\
8& \rightarrow 13
\end{cases}
\qquad {\rm L}:
\begin{cases}
\,\,9\!\!& \rightarrow 1\\
10& \rightarrow 2,3\\
11& \rightarrow 4,5\\
12& \rightarrow 6\\
13& \rightarrow 7,8.
\end{cases}
\end{align*}
The substitution $\eta$ is primitive because you `can go' from the letter 1 to any  letter and from any letter to the letter 1.
This gives irreducibility; there is primitivity because periodicity is impossible by the first rule $1\rightarrow 1,2$.

 \noindent The substitution $\eta$ has full rank because any unit vector $$e_a=(0,\dots,0,1,0,\dots,0)$$ is a linear combination of rows of the incidence matrix $M_\eta$ of $\eta$. For $a\in {\rm L}\setminus\{9\}$ this combination is trivial, and for the other letters this is exactly forced by the choice of lengths in such a way that the largest common refinement of the induced partitions of the images of S and L is the singleton partition.
In more detail:  denote the $a^{\rm th}$ row of  $M_\eta$ by $R_a$.  Then $e_1=R_9$, and thus  $e_2=R_1-R_9$, $e_3=R_{10}-e_2=R_{10}-R_1+R_9$, etc.

\smallskip

The argument yielding the property of full rank will hold in general for all $n\ge5$. To prove primitivity for all $n$ we need some more details.

\begin{proposition} The substitution $\eta=\eta_n$ is primitive for all $n \ge 5$.\end{proposition}

\noindent \emph{Proof:} The proposition will be proved if we show that for all  $a\in A$ the letter $a$ will occur in some iteration $\eta^k(1)$, and conversely,
that for all  $a\in A$ the letter $1$ will occur in some iteration $\eta^k(a)$. The first part is easy to see from the fact that $\eta(1)=1,2$ and that
$\eta^2(1,\dots,F_{n-2})=1,\dots,F_n-1$, plus $\eta^2(a_{\rm M})=F_n,1$. For the second part, we show that A) for any $a\in$ M$\cup$L a letter from S will occur in $\eta^k(a)$  in $k\le$  Card M$\cup$L steps (see Lemma~\ref{lem:dec}) and B), that  for any $a\in$ S the letter 1 will occur in $\eta^k(a)$  in $k\le$ 2Card $A$ steps (see Lemma~\ref{lem:occ1}). \hfill $\Box$

 \begin{lemma}\label{lem:dec}  Let $f:A\rightarrow A$ be the map that assigns the first letter of $\eta^2(a)$ to $a$. Then $f$ is strictly decreasing on L $\cup$ M$\backslash \{a_{\rm M}\}$.\end{lemma}

\noindent \emph{Proof:}  First we consider $f$ on ${\rm L}$. We have $$\eta^2(a_{\rm L}\dots F_n)=\eta(1,\dots, F_{n-1}-1, F_{n-1})=1\dots F_n.$$
Since $$\eta^2(F_n)=\eta( F_{n-1}-1,F_{n-1})=F_{n-1}-1+F_{n-2},F_{n-1}+F_{n-2}=F_n-1, F_n,$$
we obtain $f(F_n)=F_n-1<F_n$. This implies that also the previous letters in ${\rm L}$ are mapped by $f$ to a smaller letter.

\noindent Next we consider $f$ on M$\backslash \{a_{\rm M}\}$. Here $$\eta^2(a_{\rm M}+1,\dots, F_{n-1})=\eta(a_{\rm L}+1,\dots, F_{n})=2,3,\dots, F_{n-1}.$$
Now  $$\eta^2(F_{n-1})=\eta( F_{n})=F_{n-1}-1,F_{n-1}.$$
So we obtain $f(F_{n-1})=F_{n-1}-1<F_{n-1}$. This implies that also the previous letters in ${\rm M}$ are mapped by $f$ to a smaller letter.\hfill $\Box$

\begin{lemma}\label{lem:occ1} For all $a\in S$ there exists $k \le 2\,{\rm Card}\, A$ such that the letter 1 occurs in $\eta^{k}(a)$.
\end{lemma}

\noindent \emph{Proof:} The substitution $\eta^2$ maps intervals $I$ to intervals $\eta^2(I)$, provided $I$ does not contain $a_{\rm M}$ or $a_{\rm L}$.
By construction, since the $\eta(b)$ for $b\in {\rm L}$ have length 1 or 2, the length of $\eta(a)$ for $a\in {\rm S}$ is  2 or 3, and so $\eta(a)$ contains a word $c, c+1$ for some $c\in A$. Since $\rho\varphi^{(n-4)}(0)$ does not contain two consecutive 1's (property (V3)), the image $\eta^2(c,c+1)$ has at least length 3. Since\footnote{This follows from the fact that any word in the language of $\eta$ occurs in some concatenation of the two words $12\dots F_{n}$ and $12\dots F_{n-1}$.} any word of length at least 3 in the language of $\eta$ contains an interval of length 2,  the length increases by at least 1 if you apply $\eta^2$. It follows  that for all $n\ge 5$ and all  $a\in {\mathrm S}$ one has $|\eta^{2n+1}(a)| \ge n+2$.
But then after less than ${\rm Card}\, A$ steps a letter  $a_{\rm M}$ or  a letter $a_{\rm L}$ must occur in $\eta^{2n+1}(a)$. This implies that the letter 1 occurs in $\eta^{2n+3}(a)$, since both  $\eta^2(a_{\rm M})$ and $\eta^2(a_{\rm L})$ contain a 1.\hfill $\Box$

\section{The 2-symbol case}\label{sec:two}

The eigenvalue group of the Fibonacci system is the rotation over the small golden mean $\gamma=(\sqrt{5}-1)/2$ on the unit circle, and any system topologically isomorphic to the Fibonacci system must have an incidence matrix with Perron Frobenius eigenvalue the golden mean or a power of the golden mean
(cf.~\cite[Section 7.3.2]{Pytheas}). Thus, modulo a permutation of the symbols, on an alphabet of two symbols the incidence matrix
with Perron-Frobenius eigenvalue the golden mean has to be $\left( \begin{smallmatrix} 1 \, 1\\ 1\,  0  \end{smallmatrix} \right).$ There are two substitutions with this incidence matrix: Fibonacci $\varphi$, and reverse Fibonacci $\ravphi$, defined by
$$\ravphi: \qquad 0\rightarrow\,10,\;1\rightarrow 0.$$
These two substitutions are essentially different, as they have different standard forms (see \cite{Dekking-2016} for the definition of standard form).

However,  it follows directly from   Tan Bo's criterion in his paper \cite{Tan}
that $\ravphi$ and $\varphi$ have the same language\footnote{This follows also directly from the well-known formula  
$\ravphi^{\!2n}(0)\,10=01\,\varphi^{2n}(0)$ for all $n\ge1$ (see \cite[p.17]{Berstel}).}, 
but then they also generate the same system. Conclusion: the conjugacy class of Fibonacci with Perron-Frobenius eigenvalue the golden mean
restricted to two symbols consists of Fibonacci and reverse Fibonacci.

\section{Maximal equicontinuous factors}\label{sec:equi}

 Let $T$ be the mapping from the unit circle $Z$ to itself defined by $Tz=z+\gamma \mod 1$,
where $\gamma$ is the small golden mean.
This, being an irrational
rotation, is indeed an equicontinuous dynamical system -- the usual
distance metric is an invariant metric under the mapping.
 The factor map from the Fibonacci dynamical system $(X_\varphi,\sigma)$ to $(Z,T)$ is  given  by
requiring that the cylinder sets $\{x:x_0=0\}$ and $\{x:x_0=1\}$ are mapped to
the intervals $[0,\gamma]$ and $[\gamma,1]$ respectively, and requiring equivariance.
If we take any point of $Z$ not of the form $n\gamma \mod 1$ ($n$ any integer), then the
corresponding sequence is unique. If, however, we use an element in the
orbit of $\gamma$, then for this point there will be two codes, a ``left" one
and a ``right" one.

\smallskip

We want to understand more generally why  two or more points map to a single point. Suppose $x$ and $y$ are two
 points of a system $(X,\sigma)$ that map to two  points $x'$ and $y'$ in
an equicontinuous factor. Then for any power of $T$ (the
map of the factor system) the distance between $T^n(x')$ and $T^n(y')$ is
just equal to the distance between $x'$ and $y'$. So  $x$ and
$y$ map to the same point $x'$ if  either all $x_n$ and $y_n$ are equal for
sufficiently large $n$, or all $x_n$ and $y_n$ are equal for sufficiently
large $-n$. We say that $x$ and $y$ are respectively \emph{right asymptotic} or \emph{left asymptotic}

\smallskip

A pair of letters $(b,a)$ is called a \emph{cyclic pair} of a substitution $\alpha$ if $ba$ is an element of the language of $\alpha$, and for some integer $m$
$$\alpha^m(b)=\dots b \quad{\rm and}\quad \alpha^m(a)=a\dots. $$
Such a pair gives an infinite sequence of words $\alpha^{mk}(ba)$ in the language of $\alpha$, which---if properly centered---converge to an infinite word which is a fixed point of $\alpha^m$. With a slight abuse of notation we denote this word by $\alpha^{\infty}(b)\cdot \alpha^{\infty}(a)$.

\smallskip

For the Fibonacci substitution $\varphi$, $(0,0)$ and $(1,0)$ are cyclic pairs, and the two synchronized points $\varphi^\infty(0)\cdot\varphi^\infty(0)$ and
$\varphi^\infty(1)\cdot\varphi^\infty(0)$, are right asymptotic so they map to the same point in the equicontinuous factor.


\smallskip

Because of these considerations we now define $Z$-triples. Let $\eta$ be a primitive substitution. Call three points $x$, $y$, and $z$ in $X_\eta$ a $Z$-\emph{triple} if they are  generated by three cyclic pairs
 of the form $(b,a),\, (b,d)$ and $(c,d)$, where $a,b,c,d \in A$. Then $x$, $y$, and $z$ are mapped to the same point in the maximal  equicontinuous factor.

 \begin{theorem}\label{th:Zth} \; Let $(X_\eta,\sigma)$ be any substitution dynamical system topologically isomorphic to the Fibonacci dynamical system. Then there do not exist $Z$-triples in $X_\eta$.
\end{theorem}

\noindent \emph{Proof:} Since $(X_\eta,\sigma)$ is topologically isomorphic to $(X_\varphi,\sigma)$, its maximal equicontinuous factor is $(Z,T$), and the factor map is at most 2-to-1. Suppose $(b,a),\, (b,d)$ and $(c,d)$ gives a $Z$-triple $x,y,z$ in $X_\eta$. Noting that
$$x=\eta^\infty(b)\cdot\eta^\infty(a), \quad y=\eta^\infty(b)\cdot \eta^\infty(d)$$
are left asymptotic, and $y=\eta^\infty(b)\cdot \eta^\infty(d)$ and $z=\eta^\infty(c)\cdot\eta^\infty(d)$ are right asymptotic, this would give a contradiction. $\Box$

\medskip

\noindent {\bf Example}  Let $\eta$ be the substitution given by
$$\eta:\qquad 1\rightarrow 12,\,
 2\rightarrow 34,\,
 3\rightarrow 5 ,\,
 4\rightarrow 1,\,
 5\rightarrow 23. $$
 Then $\eta$ generates a system that is topologically isomorphic to the Fibonacci system ($\eta$ is the substitution at the end of Section~\ref{sec:reshaping}). Quite remarkably, $\eta^6$ admits 5 fixed points generated by the cyclic pairs $(1, 2),\, (2, 3),\, (3, 1),\, (4, 5)$ and $(5, 1)$.
 Note however, that no three of these form a $Z$-triple.

\section{Fibonacci matrices}\label{sec:mat}

 Let $\mathcal{F}_r$ be the set of all non-negative primitive $r\times r$ integer matrices, with Perron-Frobenius eigenvalue the golden mean $\Phi = (1+\sqrt{5})/2$.\\ 
 We have seen already that $\mathcal{F}_2$ consists of the single matrix $\left( \begin{smallmatrix} 1 \, 1\\ 1\,  0  \end{smallmatrix} \right).$

%
%


\begin{theorem}\label{th:F3} The class  $\mathcal{F}_3$ essentially consists of the three matrices

 \smallskip

\qquad $ \left( \begin{smallmatrix} 0\, 1\, 0\\ 1\, 0\, 1\\ 1\, 1\, 0 \end{smallmatrix} \right),\;
 \left( \begin{smallmatrix}  0\, 1\, 0\\  0\, 0 \,1 \\ 1\, 2\,0 \end{smallmatrix} \right),\;
 \left( \begin{smallmatrix} 0\, 1\, 0\\ 1 \,0\, 1\\ 1\, 0\, 1 \end{smallmatrix} \right).$
\end{theorem}

 Here essentially means that in each class of 6 matrices corresponding to the permutations of the $r=3$ symbols, one representing member has been chosen (actually corresponding to the smallest standard form of the substitutions having that matrix).

 \medskip

\emph{Proof:}   Let $M$ be a non-negative primitive $3\times 3$ integer matrix, with Perron-Frobenius eigenvalue the golden mean $\Phi = (1+\sqrt{5})/2$. We write\\ [-0.7cm]

 $$ M= \left( \begin{matrix} a\; b\; c\\ d\; e\; f\\ g\; h\; i \end{matrix} \right).$$
The characteristic polynomial of $M$ is $\chi_M(u)=u^3-Tu^2+Fu-D,$
where $T=a+e+i$ is the trace of $M$, and
\begin{equation}\label{eq:FandD}
F=ae+ai+ei-bd-cg-fh,\quad
D=aei+bfg+cdh-afh-bdi-ceg.\quad
\end{equation}
Of course $D$ is the determinant of $M$.
Since $\Phi$ is an eigenvalue of $M$, and we consider matrices over the integers, $u^2-u-1$ has to be a factor of $\chi_M$.
Performing the division we obtain
$$\chi_M(u)=\big(u-(T-1)\big)\big(u^2-u-1\big),$$
and requiring that the remainder vanishes, yields
 \begin{equation}\label{eq:DF}
F=T-2,\quad D=1-T.
\end{equation}
Note that the third eigenvalue equals $\lambda_3=T-1$. From the Perron-Frobenius theorem follows that this has to be smaller than $\Phi$ in absolute value, and since it is an integer, only $\lambda_3=-1, 0, 1$ are possible. Thus there are only 3 possible values for the trace of $M$: $T=0,\, T=1$ and $T=2$.

\smallskip

 The smallest row sum of $M$ has to be smaller than the PF-eigenvalue $\Phi$ (well known property of primitive non-negative matrices). Therefore $M$ has to have one of the rows $(0,0,1)$, $(0,1,0)$ or $(0,0,1)$. Also, because of primitivity of $M$, the 1 in this row can not be on the diagonal. By performing permutation conjugacies of the matrix we may then assume that $M$ has the form
 $$ M= \left( \begin{matrix} 0\;\; 1\;\; 0\\ d\;\; e\;\; f\\ g\;\; h\;\; i \end{matrix} \right).$$
 The  equation \eqref{eq:FandD}  combined with \eqref{eq:DF} then simplifies to
 \begin{equation}\label{eq:DF2}
T-2=F=ei-d-fh, \quad
1-T=D=fg-di.
\end{equation}


\noindent{\bf Case $\mathbf{{\emph T}=0}$}

 \noindent In this case $e=i=0$, so \eqref{eq:DF2}  simplifies to
 \begin{equation}\label{eq:T0F}
-2=F=-d-fh, \quad 1=D=fg.
\end{equation}
Then $f=g=1$, and so $d+h=2$. This gives three possibilities leading to the matrices
$ \left( \begin{smallmatrix} 0\, 1\, 0\\ 1\, 0\, 1\\ 1\, 1\, 0 \end{smallmatrix} \right),\;
 \left( \begin{smallmatrix}  0\, 1\, 0\\  0\, 0 \,1 \\ 1\, 2\,0 \end{smallmatrix} \right),\;
 \left( \begin{smallmatrix} 0\, 1\, 0\\ 2 \,0\, 1\\ 1\, 0\, 0 \end{smallmatrix} \right).$
%
%

\noindent Here the third matrix is permutation conjugate to the second one.

 \smallskip

\noindent{\bf Case $\mathbf{{\emph T}=1}$}

 \noindent In this case $e=1, i=0$, or $e=0, i=1$.

  \noindent First case: $e=1, i=0$. Now  \eqref{eq:DF2} simplifies to
 \begin{equation}\label{eq:T1F}
-1=F=-d-fh, \quad 0=D=fg.
\end{equation}
Then $g=0$, since $f=0$ is not possible because of primitivity.
But $g=0$ also contradicts  primitivity, as  $d+fh=1$, gives either $d=0$ or $h=0$.

 \noindent Second case: $e=0, i=1$.
Now  \eqref{eq:DF2} simplifies to
\begin{equation}\label{eq:T1F2}
-1=F=-d-fh, \quad 0=D=fg-d.
\end{equation}
Then $d=0$  would imply  that $f=h=1$.
But, as $g>0$  because of primitivity, we get a contradiction with $fg=d=0$.

 On the other hand, if $d>0$, then $d=1$ and $f=0$ or $h=0$. But $fg=d=1$ gives $f=g=1$, so $h=0$, and we obtain the matrix
 $ \left( \begin{smallmatrix} 0\, 1\, 0\\ 1\, 0\, 1\\ 1\, 0\, 1 \end{smallmatrix} \right).$


 \noindent{\bf Case $\mathbf{{\emph T}=2}$}

 \noindent  In this case  \eqref{eq:DF2} becomes
 \begin{equation}\label{eq:T1F}
0=F=ei-d-fh, \quad -1=D=fg-di.
\end{equation}
Since $ei=0$ would lead to $f=0$, which is not allowed by primitivity, what remains is $e=1,i=1$. Then, substituting $d=fg+1$ in the first equation gives
$0=f(g+h)$. But both $f=0$ and $g=h=0$ contradict primitivity.

\smallskip

Final conclusion: there are three matrices in $\mathcal{F}_3$, modulo permutation conjugacies. \hfill$\Box$

\bigskip

\noindent {\bf Remark} It is well-known that the PF-eigenvalue lies between the smallest and the largest row sum of the matrix. One might wonder whether this largest row sum is bounded for the class $\mathcal{F}=\cup_r\mathcal{F}_r$.
Actually the class $\mathcal{F}_r$ contains matrices with some row sum equal to $r-1$ for all $r\ge 3$:
  take the  matrix $M$ with  $M_{1,j}=1$ for $j=2,\dots,r$, $M_{2,2}=1$ and
 $M_{i,{i+1}}=1$, for $i=2,\dots,r-1$, $M_{r,1}=1$ and all other entries 0.


Now note that $(1, \Phi,...,\Phi)$ is a left eigenvector of $M$ with eigenvalue $\Phi$ (since $\Phi^2=1+\Phi$).
Since the eigenvector has all entries positive, it must be a PF-eigenvector (well known property of
primitive, non-negative matrices), and hence $M$ is in $\mathcal{F}_r$.

\section{The 3-symbol case}\label{sec:C2}



 \begin{theorem} There are two primitive injective substitutions
 $\eta$ and $\zeta$ on a three letter alphabet $\{a,b,c\}$  that generate dynamical systems topologically isomorphic to the Fibonacci system. These are given\footnote{Standard forms: replace $a,b,c$ by $1,2,3$.}  by\\[-.4cm]
 $$\eta(a)=b,\,\eta(b)=ca,\, \eta(c)=ba,\quad \zeta(a)=b,\,\zeta(b)=ac,\, \zeta(c)=ab.$$
\end{theorem}

\noindent\emph{Proof:}  The possible matrices for candidate substitutions are given in   Theorem~\ref{th:F3}. Let us consider the first matrix
$ \left( \begin{smallmatrix} 0\, 1\, 0\\ 1\, 0\, 1\\ 1\, 1\, 0 \end{smallmatrix} \right)$.
\noindent There are four substitutions with this matrix as incidence matrix:
\begin{align*}
\eta_1: \; a& \rightarrow b,\, b\rightarrow ca,\,c\rightarrow ba,& \eta_2: \; a \rightarrow b,\, b\rightarrow ca,\,c\rightarrow ab,\\
\eta_3: \; a& \rightarrow b,\, b\rightarrow ac,\,c\rightarrow ba,& \eta_4: \; a \rightarrow b,\, b\rightarrow ac,\,c\rightarrow ab,
\end{align*}
Here $\eta_1=\eta$. To prove that the system of $\eta$ is conjugate to the Fibonacci system consider the letter-to-letter map $\pi$ given by
$$\pi(a)=1,\quad \pi(b)=\pi(c)=0.$$ Then $\pi$ maps $X_\eta$ onto $X_\varphi$, because $\pi\eta=\varphi\pi$. Moreover,  $\pi$ is a conjugacy, since if $x\ne y$ and $\pi(x)=\pi(y)$, then there is a $k$ such that $x_k=b$ and $y_k=c$. But the words of length 2 in the language of $\eta$ are $ab, ba, bc$ and $ca$, implying that $x_{k-1}=a$ and $y_{k-1}=b$,  contradicting $\pi(x)=\pi(y)$.

 Since $\zeta$ is the time reversal of $\eta$, and we know already that the system of $\ravphi$ is conjugate  to the Fibonacci system, the system generated by $\eta_4=\zeta=\ate$ is conjugate to the Fibonacci system.

 It remains to prove that $\eta_2$ and $\eta_3$ generate systems that are \emph{not} conjugate to the Fibonacci system. Again, since $\eta_3$ is the time reversal of $\eta_2$, it suffices to do this for $\eta_2$.
 The language of $\eta_2$ contains the words $ab, bb$ and $bc$. These words generate fixed points of $\eta_2^6$ in the usual way. But these three fixed points form a $Z$-triple, implying that the system of $\eta_2$ can not be topologically isomorphic to the Fibonacci system
 (see Theorem~\ref{th:Zth}).

 The next matrix we have to consider is $\left( \begin{smallmatrix}  0\, 1\, 0\\  0\, 0 \,1 \\ 1\, 2\,0 \end{smallmatrix} \right).$
 There are three substitutions with this matrix as incidence matrix:
\begin{align*}
\eta_1: \; a& \rightarrow b,\, b\rightarrow c,\,c\rightarrow abb,& \eta_2: \; a \rightarrow b,\, b\rightarrow c,\,c\rightarrow bab,\\
\eta_3: \; a& \rightarrow b,\, b\rightarrow c,\,c\rightarrow bba.
\end{align*}
Again, the system of $\eta_1$ contains a $Z$-triple generated by $ab, bb$ and $bc$. So this system is not conjugate to the Fibonacci system, and neither is the one generated by $\eta_3$ (time reversal of $\eta_1$). The system generated by $\eta_2$ behaves similarly to the Fibonacci system, \emph{but} is has an eigenvalue $-1$ (it has a two-point factor via the projection $a,c\rightarrow 0, b\rightarrow 1$.)

Finally, we have to consider the matrix $\left( \begin{smallmatrix} 0\, 1\, 0\\ 1 \,0\, 1\\ 1\, 0\, 1 \end{smallmatrix} \right).$
 There are four substitutions with this matrix as incidence matrix:
 \begin{align*}
\eta_1: \; a& \rightarrow b,\, b\rightarrow ac,\,c\rightarrow ac,& \eta_2: \; a \rightarrow b,\, b\rightarrow ac,\,c\rightarrow ca,\\
\eta_3: \; a& \rightarrow b,\, b\rightarrow ca,\,c\rightarrow ac,& \eta_4: \; a \rightarrow b,\, b\rightarrow ca,\,c\rightarrow ca.
\end{align*}
Here $\eta_1$ and $\eta_4$ generate systems conjugate to the Fibonacci system, but the substitutions are not injective.
The substitution $\eta_2$ has all 9 words of length 2 in its language, and all of these generate fixed points of $\eta_2^6$. So the system of $\eta_2$ is certainly not topologically isomorphic to the Fibonacci system. The proof is finished, since $\eta_3$ is the time reversal of $\eta_2$.  $\Box$


\section{Letter-to-letter maps}\label{sec:L2L}

By the Curtis-Hedlund-Lyndon theorem all members in the conjugacy class of the Fibonacci system can be obtained by applying letter-to-letter maps $\pi$ to $N$-block presentations $(X^{[N]},\sigma)$. Here we analyse the case $N=2$. The 2-block presentation of the Fibonacci system is generated by (see Section~\ref{sec:Nblock}) the 2-block substitution
$$\hat{\varphi}_2(1)=12 \quad \hat{\varphi}_2(2)=3, \quad \hat{\varphi}_2(3)=12.$$
There are (modulo permutations of the symbols) three letter-to-letter maps from $\{1,2,3\}$ to $\{0,1\}$. Two of these project onto the Fibonacci system, as they are projections on the first respectively the second letter of the 2-blocks. The third is $\pi$ given by
$$\pi(1)=0,\quad \pi(2)=0, \quad \pi(3)=1.$$
What is the system $(Y,\sigma)$ with $Y=\pi\big(X^{[2]}\big)$?

First note that $(Y,\sigma)$ is conjugate to the Fibonacci system since $\pi$ is clearly invertible.
Secondly, we note that the points in $Y$ can be obtained by doubling the 0's in the points of the Fibonacci system.
This holds because $\pi(12)=00,\, \pi(3)=1$, but also $$\pi(\hat{\varphi}_2(12))=\pi(123)=001,\;\pi(\hat{\varphi}_2(3))=\pi(12)=00.$$
Thirdly, we claim that the system $(Y,\sigma)$ cannot be generated by a substitution. This follows from the fact that the sequences in $Y$ contain the word 0000, but no other fourth powers. This is implied by the $4^{\rm th}$ power free-ness of the Fibonacci word, proved in \cite{Karhumaki}.

A fourth property is that the sequence $y^+$ obtained by doubling the 0's in $\f$, where $\f$ is the infinite Fibonacci word  is given by
$$y^+_n=[(n+2)\Phi]-[n\Phi]-[2\Phi], \qquad {\rm for\;} n\ge 1,$$
according to \cite{Wolfdieter}, and \cite{OEIS-Fib} (here $[\cdot]$ denotes the floor function).

\smallskip

Finally we remark that Durand  shows in the paper \cite{Durand} that the Fibonacci system is prime \emph{modulo topological isomorphism}, and ignoring finite factors and rotation factors. This implies that all the projections are automatically invertible, if the projected system is not finite.

\bibliographystyle{plain}


\end{document}